%&amstex
\input epsf

\loadeurm

\loadeusm

\loadeufb

\magnification=\magstep1

\documentstyle { amsppt}

\def\sp{\  }
\def\slash{ /}
\def\a{ \alpha}
\def\b{ \beta}
\def\g{ \gamma}

\def\d{ \delta}
\def\D{ \Delta}

\def\f{ \theta}

\def\k{ \kappa}
\def\l{ \lambda}

\def\m{ \mu}

\def\c{ \chi}

\def\p{ \pi}

\def\si{ \sigma}
\def\Si{ \Sigma}

\def\t{ \tau}

\def\o{ \omega}

\def\8{_{ \infty}   }

\def\ex{ \text{ Ext}   }

\def\eC{ \eusm C}
\def\ee{ \eufm e}
\def\eg{ \eufm g}
\def\eT{ \eusm T}
\def\eI{ \eusm I}
\def\eE{ \eusm E}
\def\ef{ \eufm f}
\def\eu{ \eufm u}
\def\eU{ \eufm U}

\hfuzz 27pt

\topmatter \title Turaev-Viro Modules of Satellite knots  \endtitle
\rightheadtext{ Turaev-Viro Modules}    \author Patrick M. Gilmer \endauthor
 \affil Louisiana State University \endaffil
 \address Department of Mathematics,  Baton Rouge, LA
70803 U.S.A  \endaddress
  \email gilmer\@ math.lsu.edu \endemail

\abstract
Given an oriented knot $K$ in $S^3$ and a TQFT, Turaev and Viro defined modules
somewhat analogous to the Alexander module.
We work with the $(V_p,Z_p)$ theories of
Blanchet, Habegger, Masbaum and  Vogel \cite{BHMV}   for $p \ge 3,$
and consider the associated modules. In \cite{G}, we defined modules which also
depend on
the extra data of a  color $c$ which is assigned to a meridian of the knot in
the construction of the module. These modules can be used to calculate the
quantum invariants of cyclic branched covers of knots and have other uses.

Suppose now that $S$ is a satellite knot with companion $C,$ and pattern $P.$
We give formulas for the Turaev-Viro modules for $S$
in terms of the Turaev-Viro modules of $C$  and similar data coming from the
pattern $P.$ We compute these invariants explicitly in several examples.
\endabstract
\thanks This research was supported by a grant from
the Louisiana Education Quality Support Fund \endthanks
\keywords  Turaev-Viro module, quantum invariant, companion, pattern, TQFT
\endkeywords
\subjclass 57M99 \endsubjclass
  \endtopmatter

\document
\centerline{This version:(10\slash 11\slash 96); First Version: (10\slash
2\slash 96)}
%\head Introduction\endhead
\head \S1 Turaev-Viro modules\endhead
Let $(V,Z)$ be a Topological Quantum
Field Theory over a field  $f$  defined on a cobordism category
whose morphisms are oriented  $3$-manifolds perhaps
with extra structure. Let $(M,\chi)$ be a closed oriented
$3$-manifold $M$ with this extra structure together with
$\chi \in  H^1(M)$ where $\chi:H_1(M) \rightarrow \Bbb Z$ is onto.  Let $M_{
\infty}   $ denote the infinite cyclic
cover of $M$ given  by $\chi.$    Consider a fundamental
domain $E$ for  the  action of the  integers on  $M_{ \infty}   $
bounded by lifts of a surface $\Sigma$ dual to   $\chi,$ and in
general position. $E$ can be viewed as a cobordism from
$\Sigma$ to itself.  $Z(E)$ can be viewed as an endomorphism of $V(\Sigma).$

Let $\Cal K(V)$ be the generalized 0-eigenspace for the action of $Z(E)$ on
$V(\Sigma),$ i.e.  $\Cal K(V) =
\cup_{ k \ge 1}     \text{ Kernel}   ( Z^k).$ $Z(E)$ induces an automorphism
$Z^\flat  (E)$ of $V^\flat  (\Si)= V(\Sigma)\slash \Cal K(V)$. Alternatively
$\Cal
V^\flat  $ can be defined as $ \cap_{ k \ge 1}     \text{ Image}   ( Z^k).$
The Turaev-Viro module  $(M,\chi)$ associated to $(V,Z)$ is simply is
$V(\Sigma)^\flat  $ viewed as a $f[t,t^{ -1}   ]$-module
where $t$ acts by  $Z^\flat  (E).$

\proclaim{ Theorem 1.1 (Turaev-Viro)}     This module does not
depend on the choice of $E.$ \endproclaim

\demo{ Sketch of Proof}    A detailed exposition of Turaev-Viro's proof
\cite{TV}    is given in \cite{G,\S 1}. Here we give the main idea.
Suppose
$E'$ is another choice of fundamental domain. Without loss of generality we may
assume that $E'$ has been shifted by the covering transformation so that
$E$ and $E'$ are disjoint. Let $W$ denote the cobordism indicated by the
following schematic diagram for the infinite cyclic cover.
$$\epsffile{ 01.ai}   $$
As $ E \cup T(W)= W \cup E',$ we have  $Z(W) \circ Z(E) = Z(E') \circ Z(W)$
after identifying $V(\Sigma)$ with $V(T(\Sigma)),$ and $V(\Sigma')$ with
$V(T(\Sigma')).$ After dividing out  by $\Cal K(V),$ and $\Cal K(V)',$ $Z_W$
becomes
invertible and so provides a similarity between $Z^\flat  (E)$
and $Z^\flat  (E').$ \qed\enddemo

We will now specialize to the case that $(V,Z)$ is the $(V_p,Z_p)$   theory of
Blanchet, Habegger, Masbaum and  Vogel \cite{BHMV}   for $p \ge 3$.
These are combinatorial versions of the Witten-Reshetikhin-Turaev TQFTs
associated to $SU(2)$ and $SO(3)$\cite{W,RT}.  We will work over the field
of
fractions of
$k_p,$  which we denote $f_p.$ A ``color'' for this theory is a integer from
zero to $p\slash 2-2$ if $p$ is even. If $p$ is odd, a color is an even
positive integer less than or equal to $p-3.$ Here we depart from the usage in
\cite{BHMV,G}, by  assuming that colors are always even, when $p$ is odd.
One
advantage is that the tensor product axiom will always hold. We let $\eC$
denote the set of colors.

A triple of colors $\{ i,j,k\} $ is called admissible if
$i+j+k \equiv 0 \pmod{ 2}   ,$
and $i \le j+k,$ $j \le i+k,$ and $k \le i+k.$ Moreover their sum must be
small. In particular   $i+j+k \le p-4,$ if $p$ is even. Also
$i+j+k \le 2p-4,$ if $p$ is odd. Let $ \Cal A$ denote the set of admissible
triples.  Let $\Cal A(i,j)$ denote the set of colors $k$
such that $\{ i,j,k\} \in \Cal A.$  Let $\Bbb A(i)$ denote the set of ordered
pairs of colors $(j,k)$  such that $\{ i,j,k\} \in \Cal A.$ We also
let $\Cal A(i)$ denote the set  of colors $j$  such that
$\{ i,j,j\} \in \Cal A.$

The objects of our cobordism theory are oriented surfaces with colored banded
points and $p_1$-structure. A banded point is simply a point with an oriented
arc through
it. The empty set $\emptyset$ is also an object. $V(\emptyset),$ according to
the axioms of for a TQFT, is the scalar field $f_p.$ A morphism is an oriented
3-manifold with $p_1$-structure with admissibly colored trivalent banded
graphs. A
trivalent banded graph is an oriented surface which deformation retracts
to a trivalent graph. A coloring is an assignment of colors to the edges of the
core graph so that the colors at each vertex are admissible.  By closed
3-manifold, from now on we mean a morphism from $\emptyset$ to $\emptyset.$ If
$M$
is a closed 3-manifold, then $Z_p(M)$ is multiplication by a scalar which is
denoted $<M>_p.$

We let $\o$ denote
the linear combination of links in solid torus $\eta \sum_{ c \in \eC}    \D_c
e_c,$ where $\D_c$ denotes the evaluation of an unknot diagram colored $c$ and
$e_c$ denotes the closure of the Jones-Wenzl idempotent of the $c$-strand
Temperley-Lieb algebra.  $\o$ is need for the surgery formula
\cite{BHMV,\S1.C}. See the discussion at the end of the proof of
\cite{G,8.2}. If we have a linear combination
admissibly colored trivalent banded graphs in $\Bbb R^3,$  by an evaluation,
we mean the scalar obtained as in \cite{KL}    with $A$ a primitive $4p$th
root
of unity. Thus the empty graph evaluates to one. In general, if $M$ denotes the
same graph in $S^3,$ the  one point compactification of $\Bbb R^3$ with
$p_1$-structure
which extends over the 4-ball, then $<M>_p$ is $\eta$ times the evaluation of
the
graph in $\Bbb R^3.$

Let $K$ be an oriented knot in $S^3.$ Let $M(K)$ denote 0-framed surgery to
$S^3$ along $K$ equipped with a $p_1$-structure with zero $\sigma$-invariant.
Let $M(K,c)$
denote $M(K)$ with a meridian colored by a color $c.$ Let $\chi$ evaluate to
one
on
a positive meridian of $K.$ A surface dual to $\chi$ will be called a splitting
surface for $M(K).$
Let $Z_p(K,c)$ denote the Turaev-Viro module of $(M(K,c),\chi),$ thought of
as a similarity class of automorphisms for a finite dimensional vector space
over $f_p.$  Thus it may be described by a matrix $\Cal M$ with entries in
$f_p.$  In this case, we will write $Z_p(K,c)= \Cal M.$
Sometimes we may wish to describe $Z_p(K,c)$ with a matrix $\Cal M$ or
automorphism $\Cal Z$ which may
have a non trivial generalized 0-eigenspace. Then we write $Z_p(K,c)\equiv \Cal
M,$ or$\Cal Z$ as the case may be, and it is understood that $Z_p(K,c)$ is
given by
the induced map on the quotient after we divide out by this
generalized 0-eigenspace. We let $Z_p(K)$ denote $Z_p(K,0).$ It is easy to see
that  $Z_p(K,c)$ is zero if $c$ is odd, and $p$ is even. $Z_p(K,c)$ is
undefined if $c$ is odd, and $p$ is odd. It is shown in \cite{G}, that
$Z_p(K)$
is unchanged if we change the string orientation on $K.$

As motivation for studying $Z_p(K,c),$ we mention some results of \cite{G}.

\proclaim { Theorem 1.2}    If K is a fibered knot in a homology
sphere  which is a homotopy ribbon knot, then one is an eigenvalue of
$Z_p (K).$  \endproclaim

Let $M(K)_d$ denote the d-fold cyclic cover of $M(K)_d$ associated to $\chi$
with the $p_1$-structure induced from $M(K)$ by the projection.

\proclaim { Theorem 1.3}    $<M(K)_d>_p$ is the trace of $Z_p(K)^d.$ Thus
$<M(K)_d>_p$
can be computed by a linear recursion formula given by the characteristic
polynomial of $Z_p(K).$ \endproclaim

Let  $\sigma_{ \l}   (K)= \text{ Sign}    ( (1- \l) V+ (1- \bar\l) V^t),$
where $\o \in \Bbb C$ with $|\o|=1$ and $V$ is a Seifert matrix for $K.$
Following \cite{KM}, let  $ \si_d(K)= \sum_{ i=1}   ^{ d-1}    \sigma_{
\l_d^i}
(K),$
where $\l_d=e^{ 2 \p i\slash d}   .$
These are called  the total $d$-signatures of $K$. The $\si$ invariant of the
$p_1$-structure on $M(K)_d$ is $3\si_d(K).$
Let $K_d$ denote the branched cyclic d-fold cover of $S^3$ along $K$ with
a $p_1$-structure with $\sigma$-invariant $3\si_d(K).$

\proclaim{ Theorem 1.4}
$$<K_d>_p= \eta \sum_{ c \in \eC}   \D_c\sp \sp  Trace(Z_p(K,c)^k)$$
\endproclaim

Note that $Trace(Z_p(K,c)^k$ can be computed recursively from the
characteristic polynomial of $Z_p(K,c).$

Let $d_g(p,c)$ denote the dimension of $V_p$ of a surface of genus $g$ with a
single point colored $c.$ We have the following theorem of Walker's\cite{Wa1}
who proved that the rank of $Z(E)$ is an invariant of the pair $(M,\chi).$ His
work stimulated Turaev and Viro to refine his theorem and
prove (1.1). Theorem (1.5) may be used to estimate the genus of a knot.

\proclaim{ Theorem  1.5 (Walker)} If $K$ has genus $g,$
$$\text{rank} \left(Z_p(K,c)\right)\le d_g(p,c).$$
\endproclaim

Perhaps Walker did not consider the case $c \ne 0.$ I do not know.
{  \it From now on, to cut down on the
clutter of subscripts,  we will  omit the subscript $p$ from $Z$ and $< \sp \sp
>.$
}

\head \S 2 Satellite Knots\endhead

Let $C$  be an  oriented knot in $S^3,$ equipped with its standard framing.
The pushoff with this framing is a longitude which bounds in the complement.  A
pattern consists of a oriented link
of two components in $S^3.$ One component $A$ is called the axis and must be
unknotted. The other component $\eE$ is called the embellishment. Given $C$ and
a
pattern $\eE,$ a satellite knot $S$ is formed with $C$ as its companion.
Because
$C$ is framed, its tubular neighborhood comes equipped with an identification
to a standard solid torus. We give $A$ the standard framing.
The exterior of $A$ is also a standard solid torus with a knot $\eE$ in it.
$S$ is the image of $\eE$ if we replace the tubular neighborhood of $C$ by
the exterior of $A.$ More precisely, we recover the 3-sphere if we glue
the exterior of $C$ to the exterior of $A$ such that the oriented meridian of
$A$ goes to the oriented longitude of $C,$ and the oriented  longitude of $A$
goes to the oriented meridian of $C.$  Note that this gluing map is orientation
reversing. We sometimes denote $S$ by $C\star P.$

$$\epsffile{ 02.ai}   $$

The above example is the (2,1) cable of the figure eight knot.   A cable knot
is where $\eE$ is a
torus knot on the boundary of the exterior of $A.$
The winding number of the pattern is the linking number
of $A$ and $\eE.$ In the above example it is two. Since the invariant we are
calculating is actually insensitive to the string orientation of a knot, we
assume from now on that the winding number is nonnegative.

There is a long tradition in  knot theory for expressing invariants of
satellite knots in terms of invariants of the companion and the pattern. See
the following papers and the refences therein: \cite{ML}
for abelian invariants,\cite{Li}    for signatures and Casson-Gordon
Invariants,
and \cite{MS1,MS2}    for the Jones polynomial.
We mention now some precursor results from \cite{G}    on Turaev-Viro modules
of
satellite knots.

A connected sum of $K_1$ and $K_2$ can be viewed as a satellite with companion
$K_1$ and the pattern of winding number one obtained by taking $\eE$ to be
$K_2$
and $A$ to be a meridian of $K_2.$   In \cite{G,(7.4)}, we showed

$$Z(K_1\# K_2)= \bigoplus_{ c \in \eC}
Z(K_1,c)\otimes Z(K_2,c).\leqno(2.1)$$

More generally:

$$Z(K_1\# K_2,c)= \bigoplus_{ (i,j) \in \Bbb A(c)}
Z(K_1,i)\otimes Z(K_2,j).\leqno(2.2)$$

Below we will develop a formula which generalizes these.
Another important satellite construction is that of the k-twisted double.
The winding number is zero. Here is the pattern $D(k)$(with k=-1):

$$\epsffile{ 03.ai}   $$

 In \cite{G}, we derived formulas
for $Z(C \star D(k)),c).$ The formulas were rather complicated. Let $C_s$
denote the
evaluation of a diagram of $C$ with zero writhe colored by
$s.$ We will show in \S 8 how using the methods of \cite{G}, one can obtain
the following formula under the hypothesis that $C_s$ is nonzero for all colors
$s$:
$$Z (C \star D(k))\equiv \eta \k^{ -3}    \left(\mu_i \sum_{ s \in  \Cal
A(i,j)}
\mu_s^k C_s \right)_{ i,j \in \eC}   \leqno(2.3)$$
where $\mu_i= (-1)^iA^{ i^2+2i} $ is the contribution of a positive curl
colored
$i.$
Here we also use the notation that $(a_{ i,j}   )_{ i,j\in \Cal I}  $ denotes
the
square matrix with entries $a_{ i,j}  $ as $i$ and $j$ range over $\Cal I.$
More
generally, if $C_i$ is nonzero for all $i \in \Cal A(c),$ we will obtain,
$$Z (C \star D(k),c)\equiv  \eta \k^{ -3}      \left(\frac{ \D_i \mu_i}   {
\f(i,i,c)}
\sum_{ t \in
\Cal A(i,j)}
\frac{ \mu_t^k C_t}   { \f(i,j,t)}   \text{ Tet}   \bmatrix t&i&i\\c&j&j
\endbmatrix
 \right)_{ i,j \in \Cal A(c)}   . \leqno(2.4)$$
Here  we adopt some notation from \cite{KL}. $\f(i,j,k)$ denotes the
evaluation of a planar theta curve with edges colored $i,$ $j,$ and $k.$
$\text{ Tet}   \bmatrix a_1&a_2&a_3\\b_1&b_2&b_3 \endbmatrix$ denotes the
evaluation of a planar tetrahedron with edges around some face colored $a_1,$
$a_2,$ and $a_3,$ and for each $i,$  the edge opposite the edge colored $a_i$
colored $b_i.$ Note that $\text{ Tet}   \bmatrix t&i&i\\c&j&j \endbmatrix$
becomes
$\f(i,j,t)$  when we let $c=0.$ Also $\f(i,i,c)$ becomes $\D_i.$  So (2.4)
becomes (2.3).

Below we will also give a different formula for $Z (C \star D(k))$ coming from
the
satellite description. In \S 8, we also use these other formulas to give a new
derivation of (2.4) without the requirement that any $C_i$ be nonzero.

We need to give a slightly different description  of the satellite $S$ which
will be more suitable for glueing formulas. Let $m_C$ denote a meridian of $C$
in
$M(C).$  Now $m_C$ is isotopic in $M(C)$ to core of the solid torus added to
the exterior of $C$ in constructing $M(C).$ So the exterior of $m_C$ in $M(C)$
is just the exterior of $C$ in $S^3.$  But the meridian of $m_C$ is the
longitude of $C$ with the opposite orientation and the longitude
of $m_C$ using the obvious framing for $m_C$ is a meridian of $C.$ Thus we
have that:
\proclaim{ Lemma 2.5}   $S$ is the image of $\eE$ in the union of the exterior
of
$A$ with the  exterior of $m_C$ in $M(C)$ by an orientation reversing
diffeomorphism
which preserves
the longitudes but reverses the meridians. So $M(S)$ is the union of the
exterior of $A$ in $M(\eE)$ with the  exterior of $m_C$ in $M(C)$ by an
orientation reversing diffeomorphism which preserves
the longitudes but reverses the meridians.\endproclaim

\head \S 3 Gluing formulas \endhead

Let $\eT$ denote a  solid torus $S^1 \times D^2 $ with a fixed $p_1$-structure.
  The meridian $m(\eT)$ is given by $\{ 1\}
\times S^1.$ The longitude  $l(\eT)$ is given by   $S^1\times \{ 1\}  .$ Let
$\ee_i$ denote the element of $V( \eT)$ obtained by coloring
$S^1\times \{ 0\}$
with $i.$

Let $K$ a framed knot in a closed 3-manifold $M.$ By this we mean a framed knot
disjoint from the colored graph.  We may then isotope the $p_1$-structure  on
$M$
so that there is an orientation reversing diffeomorphism from a tubular
neighborhood of $K$ to $\eT $ which  preserves $p_1$-structure,
and meridian, and sends the longitude to minus the longitude. Here the
longitude is specified by the framing and the orientation on the knot. The
meridian is oriented so that $K$ intersects the meridianal disk with
intersection number plus one. We assume this has been done. Call a knot $K$
equipped with such a diffeomorphism $\phi_K,$ a pramed knot. In this situation,
let $(M,K,a)$ denote  the closed 3-manifold given by $M$ after adjoining $K$
colored $a$ to the graph. Also let $\ex (K),$ the exterior of $K$, be the
complement the interior of the domain of $\phi_K,$ with the boundary identified
with $\partial \eT$ by $\phi_K,$ which is now orientation preserving.     We
need the following lemmas.

\proclaim{ Lemma 3.1}    $\ex  (K)= \sum_{ a\in \eC}   <(M,K,a)> \ee_a$
\endproclaim
\demo{ Proof}    $ \ex (K)$ paired with $e_a$ is given by $<(M,K,a)>$ using the
pairing Q2 of \cite{BHMV}. $\{ e_a\}_{ a\in \eC}   $  is a basis and
pairing
is nondegenerate.
\qed\enddemo

\proclaim{ Lemma 3.2}   Suppose $K_1$ and $K_2$ are pramed knots in a closed
3-manifolds $M_1,$ and $M_2.$
$<\ex  (K_1)\cup_    s {
{ \phi_{ K_2}  }  ^{ -1}   \circ \phi_{ K1}
}   \left(-\ex  (K_2)\right)>= \sum_{ a\in \eC}   <(M,K_1,a)> \overline{
<(M,K_2,a)>}  $
\endproclaim
\demo{ Proof}
$<\ex  (K_1)\cup_{
{ \phi_{ { K_2}  }  }  ^{ -1}   \circ \phi_{ K_1}
}   { -\ex  (K_2)}  >= < \ex (K_1),\left(-\ex  (K_2)\right)> .$ Moreover
$\{ e_a\}_{ a\in \eC}   $  is orthonormal with respect the Hermitian pairing
Q2.
\qed\enddemo

Let $K^*$ denote $K$ in $M$ after reversing the orientation  on $M,$
but preserving the string orientation on $K.$ Note that by the above
construction starting with $K^*,$ we have $\ex(K^*)$ has boundary identified
with $\partial \eT$ by $\phi_{ K^*}  ,$ which is orientation preserving.
Also the longitude of $K^*$ is the longitude of $K,$ but the meridian of $K^*$
is oriented opposite to that of $K.$ Consider now $\ex  (K_1)\cup_    s {
{ \phi_{ K_2}  }  ^{ -1}   \circ \phi_{ K1}
}   -\ex  (K^*_2),$ which we denote by $M_1 {_{ K_1}  \wedge_{ K_2} }  M_2.$
This obtained by gluing the exterior of $K_1$ with the
exterior of $K_2$ by an orientation reversing diffeomorphism which preserves
the longitudes but reverses the meridians just as in Lemma(2.5).
(3.2) in this situation becomes:

$$< M_1 {_{ K_1}  \wedge_{ K_2} }  M_2> =\sum_{ a\in \eC}   <(M,K_1,a)>
<(M,K_2,a)>\leqno(3.3)$$

Also $M(\eE),$ and $M(C)$ both bound 4-manifolds: $B(\eE),$ and $B(C)$
respectively with zero signature such that the inclusions
$M(C) \hookrightarrow B(C),$ and $M(\eE) \hookrightarrow B(\eE)$ induce
isomorphisms on first homology. This follows from the fact that
$\Omega_3(S^1)=0.$  As
$M(\eE),$ $M(C)$ both have a $p_1$-structure with trivial $\si$-invariant, the
$p_1$-structure extends over $B(\eE),$ $B(C)$ repectively. Thus the
$p_1$-structure on
$M(\eE) {_{ A}  \wedge_{ m_C} }  M(C)$ extends over $B$ obtained by gluing
$B(\eE)$  to $B(C)$ along tubular neighborhoods of $\eE$ and $C$ identified by
the above identification.
Moreover  the kernels of the maps on $H_1$ induced by the inclusions
$\ex  (m_C)  \hookrightarrow B(C),$ and $\ex  (\eE)  \hookrightarrow B(\eE)$
are both generated by the longitudes. Thus the Maslov index of the triple of
kernels needed to compute the signature of $B$ is zero. It follows that the
$\si$-invariant of the induced $p_1$-structure on $\partial B=M(\eE) {_{ A}
\wedge_{ m_C} }  M(C)$ is zero.
Thus:

$$M(S)= M(\eE) {_{ A}  \wedge_{ m_C} }  M(C)\leqno(3.4)$$

Suppose $M_2$ above is not a closed 3-manifold, but a morphism from $\Si$ to
$\Si',$ and $K_2$ is in the interior of $M_2.$ Then one has:

$$Z( M_1 {_{ K_1}  \wedge_{ K_2} }  M_2) =\sum_{ a\in \eC}   <(M,K_1,a)>
Z(M,K_2,a))\leqno(3.5)$$

We also want to glue morphisms along the exteriors of arcs.
Let $\eI$ denote the solid tube $I \times D^2. $  We
equip $\eI$ with a fixed $p_1$-structure.  The meridian
$m(\eI)$ is given by $\{ 1\}   \times S^1,$ oriented with the standard
orientation on $S^1.$ The parallel
$p(\eI)$ is given by  $I\times \{ 1\} $  oriented with the orientation on $I$
from zero to one.

Let $M$ be a morphism from $\Si$ to $\Si',$ and suppose neither of these
surfaces are empty. Let $\g$ be a smooth framed arc in $M$ from  $p\in\Si$ to
$q \in \Si'.$  We may then isotope the $p_1$-structure  on $M$  so that there
is
an orientation reversing diffeomorphism ${ \phi_{ \g} } ^{+} $ from a
tubular neighborhood of $\g$ to
$\eI $ which  preserves $p_1$-structure, and parallel, and sends the meridian
to minus
the
meridian. Call an arc  $\g$ equipped with such a
diffeomorphism, a $(+)$pramed arc.  Here the parallel is an arc in the boundary
of the tubular
neighborhood and is specified by the framing and the orientation on the knot.
The meridian is oriented so that $\g$ intersects the meridianal disk with
intersection number plus one.

Let $-\eI$ denote $\eI$ with the same $p_1$-structure and
with the same parallel but with the opposite ambient orientation and with the
oppositely oriented meridian.
We may also isotope the $p_1$-structure  on $M$  so that there is
an orientation reversing diffeomorphism ${ \phi_{ \g} } ^{-} $ from a
tubular neighborhood of $\g$ to
$-\eI $ which  preserves $p_1$-structure, and parallel, and sends the meridian
to minus
the
meridian.   Call an arc  $\g$ equipped with such a
diffeomorphism, a $(-)$pramed arc.

In this situation, let $(M,\g,a)$
denote  morphism from say $(\Si,p,a)$ to $(\Si',q,a)$  given by $M$ after
adjoining $\g$ colored $a$ to the graph. Also let $\ex (\g),$ the exterior of
$\g$, be the complement the interior of the domain of ${ \phi_{ \g} } ^{ \pm}
,$ with the
boundary identified with $ \pm I \times S^1$  by ${ \phi_{ \g} } ^{ \pm} ,$
which is now orientation
preserving.  Suppose $\g_1$ is a $(+)$pramed arcs in a morphism $M_1,$
and  $\g_2$ is a $(-)$pramed arc in a morphism $M_2.$
 Let $M_1 {_{ \g_1} \#_{ \g_2} }  M_{ 2} $ denote $\ex  (\g_1)\cup_{
{ ({ \phi_{ \g_2} } ^+ })  ^{ -1}   \circ {( \phi_{ \g_1}^-) }
}  \left(\ex  (\g_2)\right).$ $M_1 {_{ \g_1} \#_{ \g_2} }  M_{ 2} $ is a
morphism from
$\Si_1\# \Si_2$ to $\Si'_1\# \Si'_2$ where the connect sum has been taken by
deleting neighborhoods of the endpoints of $\g_1$ and $\g_2.$ Let $\Si_{ 1,a} $
denote $\Si_{ 1} $ with the relevant point colored $a,$ etc. The  colored
splitting theorem \cite{BHMV,1.14}  describes an isomorphism $V_p(\Si_1\#
\Si_2) \approx
\oplus_{ a\in \eC}   \Si_{ 1,a} \otimes \Si_{ 2,a} .$ The following gluing
formula follows easily from the description of this isomorphism and the
definition of the maps induced by a morphism.

\proclaim{ Lemma 3.6}   If $\g_1$ is a $(+)$pramed arcs in a morphisms $M_1,$
and  $\g_2$ is a $(-)$pramed arc in a morphisms $M_2,$  then

$Z(M_1 {_{ \g_1} \#_{ \g_2} }  M_{ 2} )= \oplus_{ a\in \eC}
Z(M_1,\g_1,a)\otimes Z (M_2,\g_2,a)$
\endproclaim

Much more general gluing formulas are described in \cite{Wa2}  and the more
recent paper \cite{Ge}. The above formulas follow easily from the set-up in
\cite{BHMV}. Lemma (3.6) was used implicitly in \cite{G}  to prove (2.1) and
(2.2) above. Note also that the trace of (3.6) yields a special case of (3.3).

\head \S 4 Winding number zero  \endhead

The contribution of the companion $C$ to the formulas  we derive for $Z(S,c)$
is $<(M(C),{ m_C} ,a)>.$ This may be more convenient than the contribution
$C_a$
of $C$ to the formula (2.4) for $Z(C \star D(k),c),$ as $<(M(C),{ m_C} ,a)>=
\text{Trace}
Z(C,a).$ For instance, if $C$ itself is a satellite knot, we may use the
methods of this
paper to get our hands on$Z(C,a),$ and thus  $<(M(C),{ m_C} ,a)>.$ We will let
$C(a)$ denote
$<(M(C),{ m_C} ,a)>.$ The data $(C_a)_{ a\in \Cal C} $ and $(C(a))_{ a\in \Cal
C} $ is equivalent. In (8.1)and (8.2), we will give a change of basis matrix
which
relates these two vectors.
In this way we will rederive (2.4) from (4.4) but without
the additional hypothesis, that any $C_s$ be nonzero.

If $\eE$ has linking number  zero with the axis $A,$ then we may pick a Seifert
surface $F_P$ for $\eE$ which misses $A.$ Let $\Si_P$ denote $F_P$ capped off
in
$M(\eE).$ We may view $A$  as a subset  of $M(S),$ and $\Si_P$ misses $A.$
Let $E_{ \eE} $ be a fundamental domain for the Z-action on $M(\eE)_{
\infty}$   with boundary two lifts of $\Si,$ and continue
to call  the lift of $A$ in $E_{ \eE} ,$ by $A.$

Let $Z(\Si_P;a)=Z((E_{ \eE};a))$ where $(E_{ \eE};a)$ denotes $E_{ \eE}$
constructed as above with $A$ colored $a.$ Note that $Z^\flat((E_{ \eE};a))$
represents
the Turaev-Viro module of $M(\eE)$ with $A$ colored $a,$ which we denote
$Z(P;a).$

Similarly let $E_{ S} $ be a
fundamental domain for the Z-action on $M(S)_{ \infty} $  with boundary two
lifts of $\Si.$ Then we have
$E_{ S} = M(C){_{ m_C} \wedge_{ A}  }  E(\eE).$
So by (3.5) we have

$$Z(S) \equiv \sum_{ a\in \eC}   C(a) Z(\Si_P;a)\leqno(4.1)$$
However we cannot replace $Z(\Si_P;a)$ by $Z(P;a),$ as the summation above
requires that the maps $Z(\Si_P;a)$ have the same domain for each $a.$ Also of
course the sum of singular maps may be nonsingular.
At this point we remind the reader that when we write $Z(S)\equiv X$ where $X$
is an endomorphism or a matrix, we mean $Z(S)=X^\flat,$ where $X^\flat$ denotes
the map after dividing out by the generalized 0-eigenspace.

Similarly let $Z(\Si_P;a,c)=Z((E_{ \eE};a,c))$ where $(E_{ \eE};a,c)$ denotes
$E_{ \eE}$ constructed as above with $A$ colored $a,$ and and the inverse image
of a meridian for $\eE$ colored $c.$ Note that $Z(\Si_P;a,c))^\flat$ represents
the Turaev-Viro module of $M(\eE)$ with $A$ colored $a,$ and a meridian for
$\eE$ colored $c.$ We denote this by $Z(P;a,c).$
$$Z(S,c) \equiv \sum_{ a\in \eC}   C(a) Z(\Si_P;a,c)\leqno(4.2)$$

For the  pattern $D(k)$ of the $k$-twisted double we have:

$$Z(F_{D(k)};a)=\eta \k^{ -3}  \left(\mu_i
\sum_{ t \in \Cal A(i,j)}
{ \mu_t} ^k (-1)^{ a+t} [(a+1)(t+1)]\right)_{ i,j \in \eC}  \leqno(4.3)$$
Here $[n]$ denotes $\frac{ A^{ 2n} -A^{ -2n} } { A^{ 2} -A^{ -2} } .$ We also
have that $Z(\Si_{D(2)};a,c)$ is given by:
$$\eta \k^{ -3} \left( \mu_i \frac{ \D_i} { \f(i,i,c)}
\sum_{ t \in \Cal A(i,j)}
{ \frac { { \mu_t} ^k (-1)^{ a+t} [(a+1)(t+1)] \text{ Tet}   \bmatrix
t&i&i\\c&j&j \endbmatrix} { \f(i,j,t)} } \right)_{ i,j \in \Cal A(c)}
 \leqno(4.4)$$
Note that (4.4) becomes (4.3) if $c=0.$
Let $U$ denote the unknot. By \cite{G,7.2}, $U(a)=\d_0^a.$ Thus (4.2)
shows that $Z(U \star D(k),c)\equiv Z(P,0,c).$ Similarly $Z(U \star D(k))
\equiv Z(P,0).$ Thus

$$Z(U \star D(k))\equiv \eta \k^{ -3}  \left(\mu_i
\sum_{ t \in \Cal A(i,j)}
{ \mu_t} ^k (-1)^{t} [(t+1)]\right)_{ i,j \in \eC}  \leqno(4.5)$$

$$Z(U \star D(k),c)\equiv \eta \k^{ -3} \left( \mu_i \frac{ \D_i} { \f(i,i,c)}
\sum_{ t \in \Cal A(i,j)}
{ \frac { { \mu_t} ^k (-1)^{ t} [(t+1)] \text{ Tet}   \bmatrix
t&i&i\\c&j&j \endbmatrix} { \f(i,j,t)} } \right)_{ i,j \in \Cal A(c)}
 \leqno(4.6)$$
(4.5) becomes especially simple when $k=\pm 1.$ $U \star D(1)$ is the figure
eight knot, denoted $F8$. $U \star D(-1)$ is the right hand trefoil  knot,
denoted $RT.$

$$Z(F8)\equiv \eta \k^{ -3} \left(\mu_i^{2}\mu_j (-1)^{i+j}
[(i+1)(j+1)]\right)_{ i,j \in \eC}  \leqno(4.7)$$

$$Z(RT)\equiv \eta \k^{ -3} \left(\mu_j^{-1} (-1)^{i+j} [(i+1)(j+1)]\right)_{
i,j \in \eC}  \leqno(4.8)$$

These last two formulas, in the case p is even, are very close to the formulas
for the maps induced by the monodromies of F8 and RT \cite{G,11.1 \& 11.2}
obtained from \cite{J} after making certain substitutions and some slight
corrections. We haven't
yet seen directly that they are similar.

The rest of this section contains a derivation of (4.3), and an indication of
the proof of (4.4). We will also give the derivation of (4.7) and (4.8) from
(4.5). As in \cite{G,\S 5}, we obtain the following description
for $E_{ \Cal E} $ with the lift of $A$ colored $a$  as surgery on $S^2 \times
I$ with a tunnel drilled out from the bottom and a 1-handle added to the top.
Here and in later diagrams, the `slab' $D^2 \times I$ denotes a part of a copy
of $S^2 \times I,$ the rest of which is left out of the picture but is included
in the manifold depicted.

$$\epsffile{  04.ai} $$

Here the scalar $\k^{ -3} $ indicates that we must multiply $Z$ of the manifold
pictured by this scalar to  correct for a change in $p_1$-structure.  We may
compute the i,j
entrie of the matrix for $Z$ with respect to the basis orthonormal basis $\{
\ee_i\}_{ i \in \Cal C}  $ by adding a solid torus to bottom with core colored
$j$ and a solid torus to the top with core colored $i$ to get a colored  link
in $S^1 \times S^2.$ The invariant of this manifold is the evaluation of the
diagram on the left of Figure 4.  We expand out
the lower $\omega=\eta \sum_{ s \in \eC}    \D_s
e_s,$ and remove the $2k+1$ full twists at the cost of introducing $\mu_s^{
2k+1} .$ The sum on the right of Figure 5 is over $s \in \Cal C.$
$$\epsffile{ 05.ai} $$

Then we  use  the following simplification where the sum is over $t \in \Bbb
A(i,s).$ The second inequality uses \cite{L1,Lemma 6}   which is legitimate
when $p$ is
odd as we have chosen our colors to be even. The sum is over $t\in \Cal C.$

$$\epsffile{ 06.ai} $$

We obtain:

$$\epsffile{ 07.ai} .$$

We simplify the two strands with $k$ full twists by repeatedly using:
$$\epsffile{ 08.ai} .$$
Here the sum is over $t \in \Cal A(i,j).$ Alternatively, one may use the same
trick to take care of all k twists at once.   Then one uses, say, \cite{
KL9.10(iii)}  to collapse forks. The evaluation is seen to be
$$
\eta \k^{ -3}  { \mu_i} ^{ 2k+1}
\sum_{ t \in \Cal A(i,j)}
\left( \frac{ \mu_t} { \mu_i \mu_j} \right) ^k
H(a,t)=
{  (\frac{ \mu_i} { \mu_j} ) } ^k
 \eta \k^{ -3}  \mu_i
\sum_{ t \in \Cal A(i,j)}
{ \mu_t} ^k (-1)^{ a+t} [(a+1)(t+1)]
\leqno(4.9)$$

Here $H(a,t)$ denotes the evaluation on the standard Hopf link with components
colored $a$ and $t.$ Morton and Strickland evaluated $H(a,t)$ to be as
$(-1)^{ a+t} [(a+1)(t+1)]$ in the $p$ even case\cite{MS1,MS2}. However the same
argument works in the $p$ odd case.  There is shift by one in the index
of  corresponding colors between this paper and \cite{MS1,MS2}.  The matrix
on the right of (4.9)
is simplied by removing the factor $(\frac{ \mu_i} { \mu_j} )$ by the change of
basis $\{ \ee_j\}   \rightarrow \{ \mu_j^k
\ee_j\} .$ So we obtain the matrix given on the left hand side of (4.3).

$E_{ \Cal E} $ with the lift of $A$ colored $a$ and  the inverse image of a
meridian of $\Cal E$ colored $c$ can be pictured as in Figure 4 except one must
add a single vertical line colored $c.$
Consider the basis $\{  \ef_i\} $ for the vector space of a boundary of a solid
torus with one point colored $c$, given by the core of the solid torus colored
$i$ and an edge joining the core to the point colored $c,$ where $i \in \Cal
A(c).$ This basis is orthogonal but not orthonormal.
$Z(P;a,c)$ is given by a matrix whose $i,j$ entrie is the quotient of the
evaluations pictured in Figure 9. Here the indices $i$ and $j$ run over $\Cal
A(c).$  The same methods of evaluation then yield (4.4).
$$\epsffile{ 09.ai}  $$

To derive (4.7,) we note that $\mu_i^{-1} \mu_j^{-1} \sum_{ t \in \Cal A(i,j)}
 \mu_t  (-1)^{t} [(t+1)]= H(i,j)$ using Figure 8. Then we make use of the
Morton-Strickland formula for $H(i,j).$ (4.8) follows from the conjugate of the
above equation.

\head \S 5 Winding number one \endhead

We may  find a splitting surface $\Si_1$ in $M(C)$ which meets $m_C$ in a
single point
$x.$  Note $\Si_1$
itself may be formed from a  Seifert surface for $C$ in $S^3$ capped off in
$M(C).$ If $\eE$ has linking number one with the axis $A,$ then we may also
pick a Seifert surface $F_2$ for $\eE$ which meets $A$ in a single point say
$y.$ Let $\Si_2$ denote $F$ capped off in $M(\eE).$
$\Si_1{_{ x} \#_y}  \Si_2$ forms
a splitting surface $\Si$ in $M(S)$ assuming that the meridian around the
point $x $ is glued to the meridian around $y$ as in (2.5).
Here ${_{ x} \#_y} $
indicates that the connect summing takes places at the points $x$ and $y.$

Let $E_1$ be a fundamental domain for the Z-action on $M(C)_{ \infty} $   with
boundary two lifts of $\Si_1.$ Let $\g_1$ denote the inverse image of $m_C$ in
$E_1.$
Let $E_2$  be a fundamental domain for the Z-action on $M(\eE)_{ \infty} $
with boundary two lifts of $\Si_2.$ Let $\g_2$ denote the inverse image of $A$
in $E_2.$
Similarly let $E_{ S} $ be a fundamental domain for the Z-action on $M(S)_{
\infty} $  with boundary two lifts of $\Si.$ Then we have
$E_{ S} = E_1{_{ \g_1} \#_{ \g_2}  }  E_2.$
So by (3.6)
$$Z(S) =\bigoplus_{ a\in \eC}  Z(C,a) \otimes Z(P,a)\leqno(5.1)$$
Similarly
$$Z(S,c) =\bigoplus_{ a\in \eC}   Z(C,a) \otimes Z(P;a,c)\leqno(5.2)$$
Although there is a clear analogy between (5.1) and (5.2) and
(4.1) and (4.2), they are really quite different because of the differences
between $\sum$ and $\oplus.$ In particular the direct sum of invertible
automorphisms in still invertible but their sum need not be. Also $\oplus$ and
$\otimes$ are well behaved with respect to the operation $X \mapsto X^\flat.$

The most important example of a winding number one satellite construction is
the connect sum of two knots $K_1$ and $K_2.$ Here we take $C$ to be $K_1$
and take $E$ to be $K_2$ with the axis $A,$ a meridian to $K_2.$ Then (5.1)
yields (2.1). It is less obvious that (5.2)
yields (2.2).  Note that the meridian to $\eE$ and the axis in $P$ are
parallel.

Suppose $M$ is a morphism from a surface $\Si$ to itself which contains a
framed arc $\g$ going from point $x$ in one copy of $\Si$ to this same point
$x$ in the other copy. Suppose $\g'$ is pushoff of $\g$ going from, say, $y$ in
one copy of $\Si$ to this same point $y$ in the other copy.
Let $\Si(b)$ denotes $\Si$ with $x$  colored $b.$
Let $\Si(a,c)$ denotes $\Si$ with $x$ and $y$ colored $a$ and $c.$
Let $M(b)$ denote $M$ with $\g$ colored $b,$
$M(a,c)$ denote $M$ with $\g$ and $\g'$ colored $a$ and $c.$
We have an isomorphism from $V(\Si(a,c)$ to $ \oplus_{ b \in \Cal A(a,c)}
V(\Si(b).$ The components of this map
are given by $Z$ of a $Y$ graph colored $a,$ $b,$ and $c$ embedded in
$\Si \times I.$ Moreover it is not hard to see that under this isomorphism
$Z(M(a,c))=\oplus_{ b \in \Cal A(a,c)}  Z(M(b)).$
The above observation shows that $
Z(P,a,c)= \oplus_{ b\in \Cal A(a,c)} Z(P,b).$
Thus (5.2) implies (2.2).

We calculate now $Z(P,a)$ for certain pattern $P.$ Consider Figure 10a.
Here we have drawn $\Cal E$ as unknotted and $A$ as tangled. However one can
isotope $A$ into a standard unknot  and then $\Cal E$ becomes tangled. The
resulting picture is the pattern  $P$ we mean to study. However the link we
have
drawn is actually symmetric so, in this case, $P$ can be obtained by switching
the labels of $A$ and $\Cal E$ in Figure 10a.
$$\epsffile{ 10.ai} $$
\vskip -1.95in \hskip .75in { \bf $\Cal E$}  \vskip 1.95in

For a general pattern $\Cal E$ might not be isotopic to the unknot. Then one
would have to do surgery to $S^3$ in the complement of $P$  to unknot $\Cal E$
first. This is Rolfsen's method for calculating the Alexander module of a knot.
The resulting Figure which would play the role of Figure 10a then would have
some circles labelled  with an $\o$ and a scalar correction for the
$p_1$-structure as in
Figure 4.  Figure 10b shows $E$ where $E$ is the morphism for which
$Z^\flat  (E)= Z(P,a).$
$Z(E)$ is an endomorphism of a vector space with the basis: $\{  \ef_i\} $
indexed
by $ i \in \Cal A(a).$
In fact using methods already developed in \S 4 we have that $Z(E)f_j= \sum_i
z_{ i,j}  f_i,$ where $z_{ i,j} $ is the following quotient of two evaluations.
$$\epsffile{ 11.ai} $$

One then calculates that $z_{ i,j} $ is zero if $i \ne a,$ or if
$\{a,a,a\}\notin \Cal A.$ Note that if $\{ a,a,a\}\in \Cal A,$ then $a$ is
even.  If $\{ a,a,a\}\in \Cal A, $ we let $\g_a$
denote $z_{a,a} .$
We calculate that
$$\g_a= \frac
{  (-1)^{ \frac a 2}  A^{ \frac { -a(a+2)} { 2} }  }
{ \f(a,a,a)}
\sum_{ t \in \Cal A(a)}  \frac { \mu_t \D_t} { \f(a,a,t)}  \text{ Tet}
\bmatrix t&a&a\\a&a&a \endbmatrix$$

Let $\Cal T$ denote the set of colors $c$ such that $\{ c,c,c\} $ is
admissible.
Then $Z^\flat   (E)$ is zero if $a \notin \Cal T.$ If $a \in \Cal T,$
then  $Z^\flat   (E)$ is multiplication by $\l_a$ on a one dimensional space.
Thus by (5.1) we have:

$$Z(S)= \bigoplus_{ a \in \Cal T}  \g_a Z(C,a).$$

One could easily work out $Z(S,c)$ in a similar way, but the answer would be
more complicated.

\head \S 6 n-wheels \endhead

The type of data that a pattern with winding number greater than one
contributes
to our formula for the Turaev-Viro module of a satellite is more
complicated than an $f_p[t,t^{ -1} ]$ module or equivalently, a similarity
class of
an $f_p$-
automorphism. In this section, we define and study  $n-wheels.$ We
also construct invariants of any ordered link in $S^3$ with linking number $n$
with values in $n$-wheels of isomorphisms. The contributions of a pattern and
of
a companion to the satellite formula will be n-wheels of isomorphisms.

Let $n$ be a positive integer. An n-wheel $U$ is a sequence of  $n$ vector
space
homomorphisms
$u_i:U_i \rightarrow U_{ i+1} $ where $i\in \Bbb Z_n.$
A equivalence
from an w-wheel  $u_i:U_i \rightarrow U_{ i+1} $ to an n-wheel
$v_i:V_i \rightarrow V_{ i+1} $
is a sequence of isomorphisms
$T_i: U_i \rightarrow V_{ i+k} ,$  for some fixed $k$ such that
$u_{ i+k}  \circ T_i = T_{ i+1}  \circ w_i $ for all $i.$
If $U$ is equivalent to $V,$ we write $U \approx V.$  These conditions may
be visualized clearly as a commutative diagrams on an annulus.
Note that a 1-wheel is simply a vector space endomorphism, and equivalence of
1-wheels is similarity. The dimension of an n-wheel $u_i:U_i \rightarrow U_{
i+1} $ is simply $\dim(U_0.)$ We also have zero dimensional $n$-wheels, which
we
denote as 0. The maps of a 1-dimensional n-wheel will be denoted by scalars.

Let $\hat u_i$  denote the endomorphism of $U_i$ given by the composition $u_{
i-1}  \circ u_{ i-2}  \cdots  u_{ 0} \circ u_{ n-1}  \cdots u_{ i+1}  \circ u_{
i}  .$ Let $\Cal K(U_i)= \cup_{ k \ge 1}  Kernel(\hat {u_i}^k).$ Let ${
U_i}^\flat  = U_i \slash \Cal K(U_i).$  $u_i$ induces a map ${u_i}^\flat  $
from ${U_i}^\flat  $ to ${ U_{i+1}}^\flat  .$ In this way we get a new n-wheel
of
vector space isomorphisms.

If $U=(u_i:U_i \rightarrow U_{  i+1}  )$, and $V=(v_i:v_i \rightarrow v_{  i+1}
 )$ are two n-wheels. We may define the tensor product by
$U \otimes V=(u_i \otimes v_i:U_i \otimes V_i \rightarrow U_{  i+1}  ) \otimes
V_{  i+1}  ).$ One has the easily proved  lemma:

\proclaim {  Lemma 6.1}  $0 \otimes U =0.$
$U \otimes V \approx  V \otimes U.$
$U^\flat   \otimes V^\flat   \approx  (U \otimes V)^\flat  .$ If $U \approx U'
$ and
$V \approx V' ,$ then
$U \otimes V \approx U' \otimes V'$ \endproclaim

Given an n-wheel $U=(u_i:U_i \rightarrow U_{  i+1}  )_i, $ we may form a new
n-wheel $U(k)$ with by shifting all indices by $k,$ ie.$U(k)_i= U_{k+i} $
and $u(k)_i= u_{k+i}.$

\proclaim {  Lemma 6.2}  If $U$ is a n-wheel of isomorphisms, then
$U \approx U(k)$ for all $k.$\endproclaim

\demo{Proof} The equivalence $U \approx U(1)$ is given by the maps $u_i.$
\qed \enddemo

Given an n-wheel $U=(u_i:U_i \rightarrow U_{  i+1}  )_i, $we may form an
endomorphisms $\Bbb S(U)$ of the single vector space  $\oplus_{  i \in \Bbb
Z_n}
 U_i$ by $$\Bbb S(u_i) (\a_1 \oplus \a_2 \oplus \cdots \oplus \a_n)=
(u_n(\a_n) \oplus u_1(\a_1) \oplus \cdots u_{  n-1}  (\a_{  n-1}  ) ).$$
The similarity class of $\Cal S(U)$ is determined by the equivalence class of
$U.$

Suppose now for each $i \in \Bbb Z_n$ we have an endomorphism $g_i$ of a vector
space $G_i$, then we may form an n-wheel denoted by $\Bbb W( \{ g_i\}  ),$ as
follows.
Let $\Bbb W_i=  G_{ -i}
\otimes G_{ -i+1}  \otimes \cdots\otimes
G_{ -i+n-1} .$
Define $w_i:\Bbb W_i \rightarrow \Bbb W_{ i+1} $ by $w_i(\a_{ -i}
\otimes \a_{ -i+1}  \otimes \cdots\otimes
\a_{ -i+n-1} )=
g_{ -i+n-1} (\a_{ -i+n-1} )
\otimes \a_{ -i}  \otimes \cdots\otimes \a_{ -i+n-2} .$
This choice of indexing may seem complicated but: $\Bbb W_0=  G_0
\otimes G_1 \otimes \cdots\otimes
G_{ n-1} .$
$\Bbb W_1=  G_1
\otimes G_2 \otimes \cdots\otimes
G_{ n} ,$ and   $w_0(\a_{ 0}
\otimes \a_{ 1}  \otimes \cdots\otimes
\a_{ n-1} )=
g_{ n-1} (\a_{ n-1} )
\otimes \a_{ 0}  \otimes \cdots\otimes \a_{ n-2} .$

\proclaim {  Lemma 6.3}  $\Bbb S(U^\flat)=\Bbb S(U)^\flat.$ Also
$\Bbb W(\{ g_i^\flat\} )=(\Bbb W(\{ g_i\} ))^\flat.$
\endproclaim

Suppose we have   an oriented framed knot $K$ in a closed 3-manifold $M$ and an
epimorphism $\chi:H_1(M) \rightarrow \Bbb Z.$  Assume that $n=\chi([K])$
is greater than one.
Let $\Si$ be  surface dual to $\chi$ which meets $K$ in
exactly $n$ points.  Pick an arbitrary such point to call $x_0.$ Now travel
along
$K$ in the direction of its orientation to the next point. Call this point
$x_1.$ Continuing in this way, we may name all  $n$ points:  $\{ x_0,x_1,
\cdots, x_{ n-1}  \}  .$  Suppose we are given a ordered $w$-tuple  of colors
$\vec a
=(a_0,a_2, \cdots ,a_{ n-1} ) $ in $\Cal C^n.$ Let $F(\vec a)$ denote  $F$ with
$x_i$
colored $a_i.$  Let $\si$ denote the transformation which sends $\vec a
=(a_0,a_2,
\cdots, a_{ n-1}  ) $ to
$\si(\vec a) =(a_{ n-1} ,a_0, \cdots ,a_{ n-2}  ) .$ Let $n(\vec a)$ be the
least exponent  $e$ such that $\si^e(\vec a) = \vec a.$ Let $E$ be a
fundamental domain for the $\Bbb Z$ action on $M_\infty $  with boundary a
copies of $-\Si$ and $\Si.$ The inverse image of
$K$ consists of $n$  framed arcs.
Let $E(\vec a)$ is obtained by coloring the arc which starts at $x_i$ in $-\Si$
and goes to $x_{ i+1} $ in $F$  by $a_i,$ for all $i.$ $E(\vec a)$ is a
morphism from $F(\vec a)$ to  $F(\si \vec a).$

Let $W_p(L;\vec a)$ denote the
$n(\vec a)$-wheel given by $U^\flat  $  where $U$ is the wheel given by $U_i=
V_p(F(\si^i \vec a)),$ and $u_i=Z(E(\si^i \vec a)).$ We will usually omit the
subscript $p.$  By (6.3) $ W(L;\vec a)$
is equivalent to $ W(L,\si^i \vec a),$ so the equivalence class of $Z(L;\vec
a)$ only
depends on the cyclic ordering of $\vec a.$

  The proof of Theorem
(1.1) extends to
this situation and we have:

\proclaim{ Theorem ( 6.4)}  The equivalence class of $W(K;a)$ is an invariant
of
the isotopy class of $K.$ \endproclaim

Now suppose $L$ is a link two components $K_1$ and $K_2$ with linking number
$w$ greater than one. Let $M$ be zero framed surgery along $K_1$ with
$p_1$-structure with
zero sigma invariant, and let $K= K_2$. Define $W(L;\vec a)$ to be the $n(\vec
a)$-wheel $W(K_2;\vec a)$ defined above.
It will also be useful to do all the above with a color $c$ assigned to a
meridian of $K_2$ and the inverse images of this meridian in $E.$
In this way,
we obtain an n-wheel $W(L;\vec a,c)$ well defined up to equivalence.
We will let $U(L;\vec a,c)$ denote $U(\vec a)^\flat  $ in the above
construction.
Similarly we let $u(L;\vec a,c)$ denote the map from $U(\vec a)$ to $U(\si \vec
a).$

\head \S 7 Higher winding numbers \endhead

If view a pattern link $P$ as a link of two components where
$\Cal E$ is taken for $K_1$ and $A$ is taken for $K_2,$  and
$P$ has winding number $w,$ we
obtain $n(\vec a)$-wheels denoted $Z(P;\vec a),$ and $Z(P;\vec a,c),$ for each
$\vec a$ in $C^w.$ These only depend on $\vec a$ up to cyclic permutation.

\subhead The (2,1) cable pattern \endsubhead
We consider the
pattern $P(2,1)$ from Figure 2 with winding number two, and calculate
$W(P(2,1);\vec a,c).$     $U(P(2,1)  ,(a_1,a_2),c)$ is zero if $c \notin
\Cal A(a_1,a_2).$
If $c \in \Cal A (a_1,a_2),$  then $U(P(2,1)  ,(a_1,a_2),c)$
is one dimensional and $u(P(2,1),(a_1,a_2),c)$ is the map induced by the
manifold pictured on the left of Figure 12.
$$\epsffile{ 12.ai} $$
So $u(P(2,1)  ;(a_1,a_2),c)$ is
multiplication given by the quotient of evaluations on the left of Figure
12, which we denote by $\nu_{ a_1,a_2,c} .$  One easily has $\nu_{ a_1,a_2,c}
=\mu_{ a_1} ^{ -1}  ( \l^{ a_1,a_2}_c)^{ -1.} $   We use the $(\l^{ ab}_c)$
notation for the 3-vertex term \cite{KL, 9.9}.  Note $(\l^{ ab}_0) ^0= \d_a^b
\mu_a.$ See  \cite{MV}  for a simple derivation of $\l^{ ab}_c.$
Of course $u(P(2,1)  ;(a_2,a_1),c)$ is multiplication by
$\nu_{ a_2,a_1,c} .$

If we specialize to p=5, then $\Cal C =\{  0,2\} ,$ and
$\Cal A=\{  \{ 0,0,0\} ,  \{ 2,2,0\} ,  \{ 2,2,2\}\} .$ We have the following
non-zero
wheels for this pattern. The
one-dimensional 1-wheels $W(P(2,1);(0,0),0)=1,$ $W(P(2,1);(2,2),0)=-\bar A,$
and
$W(P(2,1);(2,2);2)= A^3.$  We also have a  one dimensional 2-wheel
$W(P;(0,2),2)$ with $u(P(2,1)  ;(0,2),2)=1,$ and $u(P(2,1) ;(2,0),2)=
A^2.$

\subhead The (3,1) cable pattern \endsubhead
Consider now the pattern $P(3,1)$ of the (3,1) cable shown on the left of
Figure 13.   This link is symmetric.  However the framed axis $A$ develops two
negative twists in the isotopy. $W(P(3,1);(a,b,c))$ is zero if $\{ a,b,c\}
\notin \Cal A.$
  If $\{ a,b,c\} \in \Cal A,$ $W(P(3,1);(a,b,c))$  is  one
dimensional. The maps (with respect to the basis given a `T' diagram with edges
labelled $a,$ $b,$ and $c$)  are multiplication by the quotient of the
evaluations on the right of Figure 13.

The numerator becomes the denominator after  removing the two kinks and two
3-vertex moves. So this quotient is $\mu_a^{ -2}  (\l^{
ab}_c)^{ -1} (\l^{ ac}_b)^{ -1} = \mu_a^{-1}.$ Thus for each
$c \in\Cal T$ we have a one wheel $W(P(1,3);(c,c,c))= \m_c^{ -1} .$ If $a \ne
b$, for each $c \in \Cal A(a,b).$ We have the 3-wheel $W(P(3,1);(a,b,c))$ with
$u(P(3,1);(a,b,c))= \mu_a^{-1}.$
$$\epsffile{ 13.ai} $$
\vskip -.9in $\quad \quad \quad \quad  \Cal E$ \vskip .9in

If we specialize to p=5, we have the following non-zero wheels for this
pattern. We have the following
one-dimensional 1-wheels or (simply vector space automorphisms).
$W(P(3,1);(0,0,0))=1,$  and $W(P(3,1);(2,2,2))=A^2.$  We also have two  one
dimensional 3-wheels:
$W(P(3,1);(0,0,2))$ with $u(P(3,1);(0,0,2))= 1,$ $u(P(3,1);(2,0,0)) =A^2,$ and
$u(P(3,1);(0,2,0))=1.$
and
$W(P(3,1);(2,2,0))$ with $u(P(3,1);(2,2,0))= A^2,$
$u(P(3,1);(0,2,2)) =1,$ and $ u(P(3,1);(2,0,2))=A^2.$

\subhead The 3-strand 1-bight turk's head  pattern \endsubhead
Consider now the pattern $T(3,1)$  shown on the left of
Figure 14.   This link is symmetric.  As the write of $\Cal E$ is zero,  no
twists develop in the framed axis $A$ during the isotopy. $W(T(3,1);(a,b,c))$
is zero if $\{ a,b,c\} \notin \Cal A.$  If $\{ a,b,c\} \in \Cal A,$
$W(T(3,1);(a,b,c))$ is  one dimensional.  The maps (with respect to the basis
given  by a `T' diagram with edges labelled $a,$ $b,$ and $c$)  are
multiplication
by the quotient of the evaluations on the right of Figure 14.

The numerator becomes the denominator after two 3-vertex moves. So this
quotient is
$(\l^{ ab}_c)(\l^{ ac}_b)^{ -1}   =\mu_b \mu_c^{-1}.$ Thus for each
$c \in\Cal T,$ we have a one wheel $W(T(1,3);(c,c,c))= 1.$ If $a \ne b$, for
each $c \in \Cal A(a,b),$ we have the 3-wheel $W(T(1,3);(a,b,c))$ with
$u(T(1,3),(a,b,c))=\mu_b \mu_c^{-1}.$
$$\epsffile{ 14.ai} $$
\vskip -.9in $\quad \quad \quad \quad  \Cal E$ \vskip .9in

If we specialize to p=5, we have the following non-zero wheels for this
pattern. We have the following
one-dimensional 1-wheels.
$W(T(3,1);(0,0,0))\allowmathbreak=1,$  and $W(T(3,1);\allowmathbreak
(2,2,2))=1.$  We also have two  one
dimensional 3-wheels:
$W(T(3,1);(0,0,2))$ with $u(T(3,1);(0,0,2))=A^2,$
$u(T(3,1);(2,0,0))=1,$
$u(T(3,1);(0,2,0))= A^8,$
and $W(T(3,1);(2,2,0))$ with
$u(T(3,1);(2,2,0))= A^8,$
  $u(T(3,1);(0,2,2)) =1,$
and $u(T(3,1);(2,0,2))= A^2.$

\subhead Wheels associated to the companion \endsubhead
Again let $\vec a=(a_0,a_2, \cdots ,a_{ w-1} )$ in $\Cal C^w.$ Then we obtain
a sequence of $n(\vec a)$ isomorphisms $Z_p(C,a_i).$ We may define the
$n(\vec a)$-wheel
of the companion associated to a sequence of  colors $\vec a$ by $\Bbb
W_p(C,\vec a)
= \Bbb W(\{ Z_p(C,a_i)\}  ).$ We will usually omit the subscript $p.$ Note that
the isomorphism class of this $n(\vec
a)$-wheel also
only depends on the cyclic ordering of $\vec a.$
We let $\eU(C,(a_1,a_2, \cdots a_{ n(\vec a)} ))$ denote the vector space
$Z(C,a_1) \otimes Z(C,a_2) \otimes \cdots
\otimes Z(C,a_{ n(\vec a)} ).$  Let  $\eu(C,(a_1,a_2, \cdots a_{ n(\vec a)}
))$ denote associated isomorphism from $\eU(C,(a_1,a_2, \cdots a_{ n(\vec a)}
))$
to $\eU(C,(a_{ n(\vec a)},a_2, \cdots a_{ n(\vec a)-1} )).$

\subhead Wheels associated to figure eight companion at p=5 \endsubhead
As an example, suppose $C$ is the figure eight knot F8.
This is the 1-twisted double of the unknot.
In \cite{G}, we calculated that has $Z_5(F8)$ two eigenvectors: $e_1$ with
eigenvalue
 $A,$ and $e_2$ with eigenvalue $\bar A.$  Also $Z_5(F8,2)$ is the identity map
on a one dimensional space. Let $ f$ denote a vector in this space.  Of
course these calculations also follow from (4.6) and (4.7).

We calculate the nonzero wheels $\Bbb W(F8, \vec a)$ associated to a color
vector $\vec a$
of length two.

So $\Bbb W(F8,(0,0))$ is a 4-dimensional 1-wheel.
$\eU(F8,0,0)$ has a basis of elements of the form
$e_i \otimes e_j $ ordered lexicographically.
With respect to this basis $\eu(F8,(0,0)$  is given by $\eufb G_1,$
the direct sum of the three matrices: $(A),$   $\pmatrix   0&A \\
\bar A & 0 \endpmatrix  $
and $(\bar A).$
%$$\eufb G_1= \pmatrix A &0&0&0\\0&0&A&0\\0&\bar A& 0&0&0\\0&0&0&\bar
% A\endpmatrix$$
 $\eufb G_1$ has eigenvalues $1,$ $-1,$ $A,$ and $\bar A.$

$W(F8,(2,2))$ is a 1-dimensional 1-wheel given by the identity.
$\eU(F8,0,2)$ has
$e_1 \otimes f, $
$e_2 \otimes f $
 as an ordered basis.
$\eU(F8,2,0)$ has
$ f \otimes e_1, $
$f \otimes e_2 $
 as an ordered basis. With respect to these bases $W(F8,(0,2))$ is a
2-dimensional 2-wheel with $\eU(F8,0,2)$ the identity and
$\eU(F8,2,0)$  given by $\eufb G_2,$ the direct sum of the two matrices: $(A),$
$(\bar A).$

We calculate now the nonzero wheels $\Bbb W(F8, \vec a)$ associated to a
color vector $\vec
a$ of length three.
$\Bbb W(F8,(0,0,0))$ is a 8-dimensional 1-wheel.
$\eU(F8,(0,0,0))$ has
$e_1 \otimes e_1 \otimes e_1,$
$e_2 \otimes e_2 \otimes e_2,$
$e_1 \otimes e_2 \otimes e_2,$
$e_2 \otimes e_1 \otimes e_2,$
$e_2 \otimes e_2 \otimes e_1,$
$e_1 \otimes e_1 \otimes e_2,$
$e_2 \otimes e_1 \otimes e_1,$
$e_1 \otimes e_2 \otimes e_1$
as ordered basis. With respect to this basis $\eu(F8,(0,0,0))$  is given by
$\eufb G_3$ the
direct sum of the four matrices: $(A),$ $(\bar A),$  $\pmatrix   0&0& A\\
\bar A&0&0 \\0&\bar A& 0 \endpmatrix , $
$\pmatrix   0&0&A \\
\bar A&0&0 \\0& A&0 \endpmatrix  .$
$\eufb G_3$ has eigenvalues  $A,$ and $\bar A,$ the three cube roots
$A,$ and the three cube roots of $\bar A.$

$\Bbb W(F8,(2,2,2))$ is a 1-dimensional 1-wheel given by the identity.

$\Bbb W(F8,(0,0,2))$ is a 4-dimensional 3-wheel.
$\eU (F8,(0,0,2))$  has a basis of elements of form
$e_i \otimes e_j \otimes f$ ordered lexicographically.
$\eU (F8,(2,0,0))$  has a basis of elements of form
$f \otimes e_i \otimes e_j $  ordered lexicographically.
$\eU (F8,(0,2,0))$  has a basis of elements of form
$e_i  \otimes f\otimes e_j $  ordered lexicographically.
With respect to these basis $\eu (F8,(0,0,2))$  is  by the
the identity matrix. $\eu (F8,(2,0,0)$ and $\eu (F8,(0,2,0)$ are both given by
$\eufb G_1.$

$\Bbb W(F8,(2,2,0))$ is a 2-dimensional 3-wheel.
$\eU (F8,(2,2,0))$  has a basis of elements of form
$f  \otimes f\otimes e_j $  ordered lexicographically.
$\eU (F8,(0,2,2))$  has a basis of elements of form
$e_i \otimes f \otimes f$ ordered lexicographically.
$\eU (F8,(2,0,2))$  has a basis of elements of form
$f \otimes e_i \otimes f $  ordered lexicographically.
With respect to these basis, $\eu (F8,(2,2,0)$   is given by $\eufb G_2.$
$\eu (F8,(0,2,2))$ and $\eu (F8,(2,0,2))$ are
given  by the
the identity matrix.

\subhead The main formula \endsubhead
Let $\Cal O^w$ denote a set consisting of a single element from  each  orbit
for
the $\Bbb Z_w$ cyclic action on $\Cal C^w,$ generated by $\si.$
We have the following generalization of (5.2). Note that it does not apply to
winding number zero patterns.

\proclaim{ Theorem (7.1)}If $S= C \star P,$ and $P$ has winding number $w,$
then:
$$  Z(S,c)= \bigoplus_{ \vec a \in   \Cal O^w} \
\Bbb S \left(  \Bbb W(C,\vec a) \otimes W(P;\vec a,c) \right)$$
\endproclaim

\demo{ Proof}  Let $\Si_C$ be a splitting surface for $M(C)$ which meets $m_C$
in a single point $x.$    Let $\Si_P$ be a splitting surface for in $M(\eE),$
which meets the axis in consecutive (along the axis) points $x_0,x_2, \cdots
x_{ w-1} .$
Let $\Si_S $ be the connected sum of $\Si_P$ at the points
$x_0,x_2, \cdots x_{ w-1} $ with $w$ copies of $\Si_C$ at the point $x.$
$\Si_P$ serves as a splitting surface for $M(S)$
in a natural way.

Let $E_C$ be the fundamental domain for the $\Bbb Z$-action on $M(C)_{ \infty}
$   with boundary two lifts of $\Si_C.$ Let $\g$ denote the inverse image of
$m_C$ in $E_C.$
Let $E_P$  be the fundamental domain for the $\Bbb Z$-action on $M(\eE)_{
\infty} $  with boundary two lifts of $\Si_P.$  The inverse image of $A$ in
$E_P$ consists of
$w$ arcs $\g_0, \g_1 \cdots \g_{ w-1} $ where $\g_i$ goes from $x_i$ in
$-\Si_P$ to
$x_{ i+1} $ in $-\Si_P.$
Let $T$ be $\Si_C \times I,$  $\t$ be the arc $\{ x\}  \times I$ in $T.$
Let $E_{ S} $ be the fundamental domain for the $\Bbb Z$-action on $M(S)_{
\infty} $  with boundary two lifts of $\Si_P.$ Then we have:

$$E_{ S} =(\cdots(({E_P}_{ \g_{ w-1} }  \wedge_{ \g} E_C)_{ \g_0} \wedge_{ \t}
T)_{ \g_1}  \wedge \cdots  )_{ \g_{ w-2} } \wedge_{ \t}  T $$

The result then follows from (3.6)
\qed \enddemo

\subhead The (2,1) cable of the figure eight \endsubhead
We consider first the case p=5. $\Cal O$ has three elements $(0,0),$ $(2,2)$
and $(2,0).$ We have $Z_5(S)$ is the direct sum of contributions from each
element of $\Cal O.$  $(2,0)$ contributes zero, as
$W(P(2,1);(0,2)=0.$  $(2,2)$ contributes $\Bbb W(F8,(2,2)) \otimes
W(P(2,1);(2,2)).$  This is a one dimensional vector space with eigenvalue
$-\bar
A.$ $(0,0)$ contributes $\Bbb W(F8,(0,0)) \otimes W(P(2,1);(0,0)) .$  This is a
4-
dimensional vector space with eigenvalues: $1$ $-1,$ $A,$  and $\bar A.$
So the eigenvalues of $Z_5(S)$ are  $1,$ $-1,$ $A,$  $\bar A,$ and $-\bar A.$

Consider now just the contribution of $\vec 0=(0,0)$ to $Z_p(S)$
$W_p(P(2,1); \vec 0)$ is the identity map on a one dimensional space. So the
contribution is just the 1-wheel $\Bbb W_p(F8,(0,0)).$ In general, we know that
$Z(F8)$ is  a unitary matrix, and so is diagonalizable with
eigenvalues all of norm one as F8 is fibered  and has genus one.
Since $F8$ is amphichiral, the non-real eigenvalues come in conjugate pairs.
Consider a pair $e_1$ and $e_2$ of eigenvectors with  conjugate
and therefore inverse eigenvalues. Then the automorphism restricted to the
subspace spanned by
$e_1 \otimes e_2$ and $e_2 \otimes e_1$ is direct summand with eigenvalues
$1$ and $-1.$ Thus in general $Z_p(F8\star P(2,1))$ has one and minus one among
its eigenvalues.

In his thesis,
Miyazaki showed $F8\star P(2,1)$ was not a ribbon knot
\cite{M}. It is  an algebraically slice fibered knot. He showed that the
monodromy does not extend over any handlebody. If the knot were ribbon, a
theorem of Casson and Gordon asserts that the closed off monodromy of a fibered
homotopy ribbon knot must extend over a handlebody \cite{CG}.  This same
theorem was used to prove Theorem(1.2). We had hoped to recover Miyazaki's
result using Theorem (1.2). So far we have not been able to do this. We working
on some refinements of Theorem(1.2). Perhaps using another TQFT would also
help.

\subhead Two winding number three satellites of the figure eight knot at
p=5\endsubhead These examples will illustrate the operation $\Bbb S$ in a
nontrivial way.
We consider first $Z_5(F8 \star P(3,1)).$
The contribution of $(0,0,0)$ is just the 1-wheel $\Bbb W(F8,(2,2,2)),$ given
by $\eufb G_3.$
The contribution of $(2,2,2)$ is just $W(P(3,1),(2,2,2)),$ given by $(A^2).$
The contribution of orbit of $(0,0,2)$ is given by the block matrix
$$\pmatrix
 0&0 & \eufb G_1\\
I_{4\times 4}&0&0\\
0&A^2 \eufb G_1 &0 \endpmatrix $$
The contribution of orbit of $(2,2,0)$ is given by the block matrix
$$\pmatrix
 0&0 & A^2 I_{2\times 2}\\
A^2  \eufb G_2&0&0 \\
0&I_{2\times 2} &0 \endpmatrix .$$
In particular, the characteristic polynomial of $Z_3( F8 \star P(3,1))$
is:
$\left(x-A\right)
\left(x-\bar A\right)
\allowmathbreak
\left(x^3-A\right)
\allowmathbreak
\left(x^3-\bar A\right)
\left(x-A^2\right)
\left(x^6 + (1- A^3)x^3-A^3\right)
\mathbreak
\left(x^{12} - (A^3 +  A^2  + A) x^9  + 2 (A^3-1)x^6  + (A^2 + A + 1)x^3
 -  A^3\right)$

$Z_5(F8 \star T(3,1))$ way be worked out in the same way as
$Z_5(F8 \star P(3,1)).$ We will just give its characteristic polynomial:
$\left(x-A\right)
\left(x-\bar A\right)
\left(x^3-A\right)
\left(x^3-\bar A\right)
\left(x-1\right)
\allowmathbreak
\left(x^{12} -(A^3 + A^2 + A)x^9 + 2(A^3 -1) x^6 + (A^2+A+1) x^3-A^3\right)
 \mathbreak  \left(x^6 + (A^3-A^2-1) x^3+ 1\right)$

\head \S 8 Derivations of (2.4) \endhead

\subhead Derivation along the lines of \cite{G} \endsubhead
  We need the hypothesis
that $C_s$ are nonzero for  $s \in \Cal A(c)$
for this approach.
Consider the exterior of the loop labelled $i$ in Figure 8a of that paper.
Instead of completing it to a diagram in $S^3,$ we should instead complete it
to a diagram in $S^1 \times S^2.$ In this way one may make use of the
orthogonality of the bases described in \cite{BHMV,4.11}. The matrix one
then
needs to invert is then already diagonal. In fact we have the following
variant of Theorem (7.7) of \cite{G}. $Z(C \star D(k),c)$ is given by matrix
whose
i,j entrie is quotient of the evaluations in Figure 15 below. Here
$i,j$ range over $\Cal A(c).$

$$\epsffile{ 15.ai}   $$

Let $\eurm B(c)_{ i,j}   $ be the evaluation of the numerator and
$\eurm L(c)_i$ be the evaluation of the denominator.
The grey disk labelled $C$ indicates that a string diagram for $C$ with zero
writhe should be inserted. The grey disk labelled $C(k)$ stands for two
parallel strands along a sting diagram for $C$
with zero writhe  with $k$ full additional twists between the strands added.
$\kappa^{ -3}   $ simply means multiply that scalar times the evaluation of the
rest of diagram.
To evaluate expand the lower loop labelled $\o$ using $\o=\eta \sum_{ s \in
\eC}
  \D_s e_s.$
Then we  use  the simplification of Figure 6.

 Now one has a strand colored $i$
with $2k+1$ full twists and $C$ with zero writhe tied into it. The twists
contribute a factor of $\mu_i^{ 2k+1}   ,$ and $C$ contributes $\frac{ C_i}
{ \D_i}   .$  Thus $\eurm B(c)_{ i,j}   = \k^{ -3}   \mu_i^{ 2k+1}   \frac{
C_i}   { \D_i}   $
times the evaluation of:

$$\epsffile{ 16.ai}   $$

Now we could remove the $k$ full twists embedded in $C(k)$ one at a time using
Figure 8.  However we can use the same trick
moving all of $C$ onto the strand labelled $t$
to show:

$$\epsffile{ 17.ai}   $$

Thus $$\eurm B(c)_{ i,j}  = \k^{ -3}   \mu_i^{ 2k+1}   \frac{ C_i}  { \D_i}
\sum_{ t \in \Bbb
A(i,j)}
\frac{ C_t}  { \theta(i,j,t)}  (\frac{ \mu_t}  { \mu_i \mu_j}  )^k \text{ Tet}
\bmatrix
t&i&i \\ c&j&j \endbmatrix.$$
Similarly
$$\eurm L(c)_{ i}  = \frac{  C_i \theta(i,i,c)}  { \eta \D_i^2}  .$$
Note $\eurm L(c)_{ i} $ is invertible if and only iff $C_i$ is invertible.
Thus $\eurm L(c)^{ -1}  \eurm B(c)$ is given by $(\frac{ \mu_i}  { \mu_j}  )^k$
times  the left
hand side of (2.4). Under the change of basis $\{ \ef_j\}   \rightarrow \{
\mu_j^k
\ef_j\}  ,$
we obtain the matrix given on the left hand side of (2.4).

\subhead A change of basis for the vector space of a torus \endsubhead
Let $T^2$ be the boundary of a solid torus. Above we have made use of the basis
 $\{\ee_i\}_{i \in \Cal C}$ for $V(T^2)$ where $\ee_i$  given by coloring a
framed core for the solid torus $i.$
Consider the basis  $\{\eg_i\}_{i \in \Cal C}$ for $V(T^2)$ where $\eg_i$ is
given by labelling a framed core for the solid torus with $\o$ and coloring the
meridian $j.$ Then $<\eg_i,\ee_j>_{T^2}$ is given by the evaluation of Figure
18.
Using Figure 6, this evaluates to $\eta H(i,j)=\eta (-1)^{i+j}[(i+1)(j+1)].$
Thus $\eg_i= \eta \sum_{j \in \Cal C}(-1)^{i+j}[(i+1)(j+1)]\ee_j.$
Applying this in a tubular neighborhood of a knot $C,$ and recalling that
$\eta^2= \frac{-(A^2-A^{-2})^2} {p},$ we obtain

$$C(i)=
\frac{-(A^2-A^{-2})^2} {p}
\sum_{j \in \Cal C}(-1)^{i+j}[(i+1)(j+1)]C_j \leqno(8.1)$$
The extra factor of $\eta$ comes from the fact that the invariant of $S^3$ with
standard $p_1$-structure
with $C$ colored $j$ is $\eta C_j.$ Now it is shown that
$\eta ((-1)^{i+j}[(i+1)(j+1)])_{i,j \in \Cal C}$ is equal to its own inverse
in \cite{MS2} in the case that $p$ is even. One can check that this is true in
general. Actually the symmetry of the above picture shows that $\ee_j= \eta
\sum_{i \in \Cal C}(-1)^{i+j}[(i+1)(j+1)]\eg_i.$ This shows that
$\eta ((-1)^{i+j}[(i+1)(j+1)])_{i,j \in \Cal C}$ is its own inverse!
Inverting the above equation, then yields
$$C_t=  \sum_{a \in \Cal C}(-1)^{a+t}[(a+1)(t+1)]C(a).\leqno(8.2)$$

$$\epsffile{ 18.ai}   $$

\subhead Derivation of (2.4) from (4.2)\& (4.4) \endsubhead
Substitute (4.4) into (4.2). Interchange the order of summation.
Then making use of (8.2), we obtain (2.4) without any hypothesis on
$C_t.$

\enddocument
\vfill\eject \end